\documentclass[twoside]{IEEEtran}
\usepackage{atbegshi,picture,ifpdf,setspace}
\usepackage[noadjust]{cite}
\usepackage[pdftex]{graphicx}

\usepackage[cmex10]{amsmath}
\usepackage[caption=false,font=footnotesize]{subfig}
\usepackage{url}

\usepackage[english]{babel}
\usepackage[utf8x]{inputenc}
\usepackage[T1]{fontenc}
\usepackage[scaled=0.91]{helvet}
\usepackage[cmintegrals,varvw,]{newtxmath}
\usepackage{acronym}
\usepackage{bm}
\usepackage{xcolor}
\DeclareMathAlphabet{\mathbfsfit}{\encodingdefault}{\sfdefault}{bx}{it}

\renewcommand*\vec[1]{\small \textsf{\textit{\textbf#1}}}

\IEEEeqnarraydefcolsep{0}{\leftmargini}

\makeatletter
\def\ps@IEEEtitlepagestyle{%
  \def\@oddfoot{\mycopyrightnotice}%
  \def\@oddhead{\hbox{}\@IEEEheaderstyle\leftmark\hfil\thepage}\relax
  \def\@evenhead{\@IEEEheaderstyle\thepage\hfil\leftmark\hbox{}}\relax
  \def\@evenfoot{}%
}
\def\mycopyrightnotice{%
  \begin{minipage}{\textwidth}
  \centering \scriptsize
  Copyright~\copyright~2018 IEEE. Personal use of this material is permitted. Permission from IEEE must be obtained for all other uses, in any current or future media, including\\reprinting/republishing this material for advertising or promotional purposes, creating new collective works, for resale or redistribution to servers or lists, or reuse of any copyrighted component of this work in other works by sending a request to pubs-permissions@ieee.org.
  \end{minipage}
}
\makeatother

\AtBeginShipout{\AtBeginShipoutUpperLeft{%
  \put(\dimexpr\paperwidth-1cm\relax,-.55cm){\makebox[0pt][r]{
  \begin{minipage}{\textwidth}
  \centering \scriptsize 
  This is the author's version of an article that has been published in this journal. 
  Changes were made to this version by the publisher prior to publication. 
  \\[0.3ex]
  The final version of record is available at\qquad http://dx.doi.org/10.1109/LAWP.2018.2869926
  \end{minipage}
  }}}%
}

\begin{document}
\title{A Weak-Form  Combined Source Integral Equation with Explicit Inversion of the Combined-Source Condition}
\author{Jonas Kornprobst, \IEEEmembership{Student Member, IEEE}, and Thomas F. Eibert, \IEEEmembership{Senior Member, IEEE}%
\thanks{Manuscript received February 2, 2018; revised August 29, 2018; accepted August 29, 2018; date of this version August 29, 2018. \emph{(Corresponding author: Jonas Kornprobst.)}}%
\thanks{The authors are with the Chair of High-Frequency Engineering, Department of Electrical and Computer
Engineering, Technical University of Munich, 80290 Munich, Germany (e-mail: j.kornprobst@tum.de).}%
\thanks{Color versions of one or more of the figures in this letter are available online
at http://ieeexplore.ieee.org.}%
\thanks{Digital Object Identifier 10.1109/LAWP.2018.2869926}%
}

\markboth{IEEE ANTENNAS AND WIRELESS PROPAGATION LETTERS}
{Kornprobst and Eibert: Explicit Inversion of the Weak-Form CS Condition}

\maketitle

\begin{abstract}
The combined source integral equation (CSIE) for the electric field on the surface of a perfect electrically conducting   scatterer can be discretized very accurately with lowest-order Rao-Wilton-Glisson  basis and testing functions if the combined-source (CS) condition is enforced in weak form. 
We introduce a technique to accelerate the iterative solution for this kind of CSIE. 
It is demonstrated that the iterative solution of the equation system can be performed very efficiently by explicitly inverting the weak-form CS condition in any evaluation  of the forward operator. 
This reduces the number of unknowns and results in improved convergence behavior at negligible, linear cost. 
Numerical results demonstrate that the new CSIE outperforms the classical CFIE for high-accuracy simulations.
\end{abstract}


\section{Introduction} 

\IEEEPARstart{E}{lectromagnetic} radiation and scattering problems involving perfect electrically conducting (PEC) objects are often treated by boundary integral equations (IEs) due to their low discretization effort and good accuracy. 
Common formulations for the tangential field components on the surface of the scatterer are the electric field IE (EFIE) and the magnetic field IE (MFIE)~\cite{Harrington1989}. 
To obtain a system of equations with one unique solution, i.e.\ avoiding the interior resonance problem, a combination of the two equations into the combined field IE (CFIE)~\cite{Mautz1979h} is the most common solution.

However, a combination of electric and magnetic sources to the combined source IE (CSIE)~\cite{Mautz1979} is also possible; 
this involves the EFIE and MFIE operators as well and is sometimes known under a different name, the Brackhage-Werner trick~\cite{brakhage1965dirichletsche,buffa2005regularized,darbas2006generalized,steinbach2009modified,melenk2012mapping}. 
A related concept (without physical currents) are the single-source surface IEs (SSSIEs), which commonly only employ one type of current unknowns~\cite{Marx1982,Glisson1984,Yeung1999,Menshov2013,Shi2015,Patel2017,Lori2018}. 

With lowest-order div-conforming expansion functions on a triangular mesh, the Rao-Wilton-Glisson (RWG) functions~\cite{Rao1990}, the discretization of the MFIE operator and also of the CFIE, is not accurate~\cite{Ergul_2004}. 
Therefore, the Buffa-Chris\-tian\-sen (BC) functions, which reside in the dual space of the RWG functions, have been successfully employed as testing functions for an accurate mixed discretization of the MFIE~\cite{Cools_2011}. 
For the CSIE as well, BC basis functions have been demonstrated to produce accurate results~\cite{Pasi2012}. 
Rather recently, we demonstrated that an accurate and interior-resonance-free formulation of the CSIE with pure RWG discretization is also possible with an RWG-tested EFIE~\cite{Kornp2017}. 
Thereby, the tangential electric field $\bm n \times\bm n \times \bm E$ is tested with RWG functions and the combined source~(CS) side condition is also implemented with RWG and $\bm n\times$RWG functions. 
Both electric and magnetic equivalent surface current densities are modeled with the div-conforming RWG basis functions, which leads to double the number of unknowns.

In this letter, we employ the CSIE discretization first in a dual-source manner, but reduce it then to a discrete SSSIE: 
the weak-form CS condition is explicitly solved in an iterative solver, leading to faster iterative solver convergence. 
In addition, the number of unknowns is reduced by a factor of two. 
First, we present the formulation and discretization of the CSIE. 
Next, the discretization and implementation strategy of the CS condition is discussed. 
Following that, simulation results demonstrate the excellent accuracy and improved iterative solver convergence and conditioning of the novel formulation of the RWG-discretized CSIE. 

\section{The Combined Source Integral Equation} 
\subsection{Formulation}
According to the Huygens principle, a PEC scattering scenario is described by the boundary condition
\begin{equation}
\bm n(\bm r) \times \bm E(\bm r) =  \bm n(\bm r) \times \bm E^\mathrm{inc}(\bm r)  + \bm n(\bm r) \times \bm E^\mathrm{sca}(\bm r)  =\mathbf{0}\,,\label{eq:PEC-Bound}
\end{equation} 
i.e.\ the superposition of the incident field $\bm E^\mathrm{inc}$ and the scattered field  $\bm E^\mathrm{sca}$ vanishes on the surface $S_0$ of the scatterer. 
The CSIE employs both  electric and magnetic surface current densities $\bm J_\mathrm S(\bm r) $ and $\bm M_\mathrm S(\bm r) $ on $S_0$ as equivalent sources to model the scattered fields.
Expressing  $\bm E^\mathrm{sca}$ in terms of the equivalent currents for a suppressed time factor $\mathrm e^{\,\mathrm j \omega t}$, the EFIE
\begin{IEEEeqnarray}{R}
\bm n(\bm r) 
\times  \bm E ^\mathrm{inc}(\bm r) 
=
\bm n(\bm r) 
\times 
\oiint\nolimits_{S_0}\bm\nabla G_0(\bm r,\bm r')
\times\bm M_\mathrm{S}(\bm r')
\mathop{}\mathrm{d}s'
\hspace*{1.2cm}
       \nonumber\\
+\frac{1}{2}\bm M_\mathrm{S}(\bm r)
       +
\mathrm{j} k_0 Z_0 \bm n(\bm r) 
\times
\left[\frac{1}{k_0^{2}}\bm\nabla
\oiint\nolimits_{S_0}G_0(\bm r,\bm r')
\bm\nabla'\cdot\bm J_\mathrm{S}(\bm r')
\mathop{}\mathrm{d}s'
          \right.\hspace*{-0.1cm}\nonumber
          \\\left.
         +
\oiint\nolimits_{S_0}G_0(\bm r,\bm r') 
\bm J_\mathrm{S}(\bm r')
\mathop{}\mathrm{d}s'
\right]\IEEEeqnarraynumspace
\label{eq:CS-EFIE}
\end{IEEEeqnarray}
is a more detailed version of~\eqref{eq:PEC-Bound}. 
The involved scalar Green's function of free space is
\begin{equation}
G_0(\bm r,\bm r')=\frac{\mathrm{e}^{-\mathrm jk_0\left|\bm r -\bm r'\right|}}{4\uppi\left|\bm r -\bm r'\right|}
\end{equation}%
with source coordinate $\bm r'$ and observation coordinate $\bm r$, 
the free-space wave number $k_0$ and the free-space wave impedance~$Z_0$. 
Since the equivalent sources are not the Love currents, and thus ambiguous, the CS constraint
\begin{equation}
\bm M_\mathrm{S}(\bm r)=\alpha Z_0\bm n(\bm r) \times \bm J_\mathrm{S}(\bm r)\label{eq:CScond}
\end{equation}
can be enforced to obtain a unique combination of electric and magnetic surface current densities, i.e.\ a quadratic system of equations with the same number of unknowns as equations. 
The combination parameter $\alpha$ is typically chosen as unity to achieve balanced contributions of electric and magnetic currents to the radiated fields~\cite{Mautz1979}. 

\subsection{Discretization}
The discretization of~\eqref{eq:CS-EFIE}  with RWG functions $\bm \beta_n$, defined on the $n$th pair of neighboring triangles, is well-known. 
Suitable testing functions in the dual space of the electric field $\bm n \times \bm n \times \bm E$ are $\bm\beta$ functions. 
The MoM equation system for PEC objects is written as
\begin{equation}
\left[-\frac{1}{2} \mathbfsfit A +  \mathbfsfit K\right]\vec v +
\mathrm j k_0 Z_0  \mathbfsfit T\,\vec i = \vec e
\end{equation}
with the matrix elements of $\mathbfsfit K$ and $\mathbfsfit T$ defined in~\cite{Ismat2009MoM} and the currents discretized as a summation of RWG basis functions 
\begin{equation}
\bm J _ \mathrm{S} = \sum\nolimits_{n=1}^N\bm \beta_n \textit{\textsf{\small i}}_n\,,
\qquad
\bm M _ \mathrm{S} = \sum\nolimits_{n=1}^N\bm \beta_n \textit{\textsf{\small v}}_n\,.\label{eq:currents}
\end{equation}
	
Furthermore, we discretize the combined source condition~\eqref{eq:CScond} in the same way, by expanding the currents according to~\eqref{eq:currents} and testing the equation with RWG functions. 
This gives the additional equation system
\begin{equation}
\mathbfsfit {A'}\vec v -\alpha Z_0\mathbfsfit A\vec i = \vec 0,\label{eq:CSIE-side}
\end{equation}
where the Gram matrix $\mathbfsfit{A}$ linking $\bm \beta$ and $\bm n\times\bm\beta$ functions has the entries
\begin{equation}
 \textit{\textsf{A}}_{mn}  =  \iint\nolimits_{S_m}\!\bm\beta_m(\bm r)\cdot\big(\bm n\times\bm \beta_n(\bm r)\big)\mathop{}\!\mathrm{d}s\label{eq:matA}
\end{equation}
 and the Gram matrix $\mathbfsfit{A}'$ of the RWG functions is defined with the entries
\begin{equation}
 \textit{\textsf{A}}_{mn}'  =  \iint\nolimits_{S_m}\!\bm \beta_m(\bm r)\cdot\bm\beta_n(\bm r)\mathop{}\!\mathrm{d}s.\label{eq:matAp}
\end{equation}
Overall, we obtain the complete CSIE equation system
\begin{equation}
\begin{bmatrix}
 \mathrm j k_0 Z_0 \mathbfsfit T&
-\frac{1}{2} \mathbfsfit A +  \mathbfsfit K \\
-\alpha Z_0 \mathbfsfit A & \mathbfsfit {A'}
\end{bmatrix}
\begin{bmatrix}
 \vec i\\
 \vec v
\end{bmatrix}
=
\begin{bmatrix}
 \vec e\\
 \vec 0
\end{bmatrix}
.\label{eq:CSIE-saddel}
\end{equation}
Since~\eqref{eq:CSIE-saddel} has electric and magnetic current unknowns, it is referred to as CSIE-JM in the following. 
In this formulation, as it has been presented in our previous work~\cite{Kornp2017}, the discretized weak-form CS condition~\eqref{eq:CSIE-side} augments the EFIE as a side condition. 
As it will be demonstrated by numerical simulations, this formulation has conditioning disadvantages compared to a formulation with only one of the sources as unknowns. 
It is important to note that the conditioning of~\eqref{eq:CSIE-saddel} is improved considerably if the magnetic current unknowns are replaced by the weighted version
\begin{equation}
\vec v' =  \vec v / Z_0 \label{eq:block}
\end {equation} 
and the correct magnetic current coefficients are retrieved in the end by a multiplication with $Z_0$. 
Alternatively, a diagonal preconditioner shows similar benefits.

To improve the implementation of the CSIE, the Gram matrices are stored separately and only the symmetric matrices~$\mathbfsfit K$ and~$\mathbfsfit T$ remain to be stored for the system matrix. 
The memory consumption of the Gram matrices is negligible with linear complexity, since $\mathbfsfit{A}$ contains $4N$ and $\mathbfsfit{A}'$ $5N$ real-valued entries. 
To reduce the number of unknowns, and halve the memory consumption of the iterative solver of the whole system, which is a generalized minimum residual (GMRES) solver, the side condition~\eqref{eq:CSIE-side} is solved explicitly within each matrix vector product according to
\begin{equation}
 \vec v =\alpha Z_0 \mathbfsfit {A'}^{-1} \!\mathbfsfit A\mathbfsfit i\,,\label{eq:pseudo}
\end{equation}
where ${\mathbfsfit {A}'}^{-1}$ denotes the inverse of~${\mathbfsfit {A}'}$ which is also obtained in an iterative manner by using the conjugate gradient method. 
The Gram matrix $\mathbfsfit {A'}$ is well-conditioned and the solver (with diagonal preconditioning) needs only a few iterations for a low residual error. Due to the linear matrix size, the convergence time of this inner iterative solver is negligible.  
We obtain the CSIE with electric current unknowns only (CSIE-J)
\begin{equation}
\left[\alpha Z_0\left(-\frac{1}{2} \mathbfsfit A +  \mathbfsfit K\right)\mathbfsfit {A'}^{-1} \mathbfsfit A  
+
\mathrm j k_0 Z_0 \mathbfsfit T\right]\vec i 
= 
\vec e\label{eq:csiej}
\,.
\end{equation}
The magnetic currents are calculated in a post-processing step according to~\eqref{eq:pseudo} since it is not necessary to store them during the solution process.

In general, the easiest way to accelerate the iterative solver convergence is diagonal preconditioning. 
For the CSIE-J, we investigate two possibilities. 
One possibility is to neglect the influence of the magnetic currents on the main diagonal and use only the diagonal matrix 
\begin{equation}
\mathbfsfit D =\mathrm j k_0 Z_0\,\mathrm{diag}(\mathbfsfit{T})
\end{equation}
for the preconditioning of~\eqref{eq:csiej}. 
Since the inverse of $\mathbfsfit{A}'$ is not directly available, only the main diagonal is taken to estimate the inverse. Hence, the matrix
\begin{equation}
\mathbfsfit D =\mathrm j k_0 Z_0 \mathrm{diag}(\mathbfsfit{T}) - \alpha Z_0\mathrm{diag}\left(\mathbfsfit{A}\,\mathrm{diag}(\mathrm{diag}(\mathbfsfit{A}'))^ {-1}\mathbfsfit{A}\right)\label{eq:CSIE-diag}
\end{equation}
can alternatively be employed for preconditioning.

\section{Numerical Results} %
For evaluating the performance of the two forms of the CSIE, they are compared to the accurate but ill-conditioned EFIE and to the well-conditioned but inaccurate MFIE. 
The common CFIE with combination parameter $0.5$ is also considered, which is comparable to $\alpha=1$ in the CSIE. 

\begin{figure}[tp]
\centering
\includegraphics{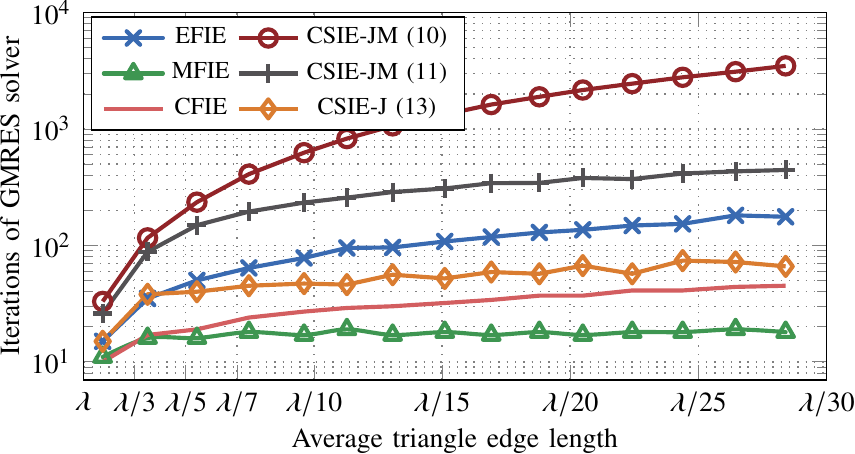}
\caption{Iterative solver convergence analysis for the mesh refinement of a 1\,m square cube.\label{fig:c-conv}}
\end{figure}

\begin{figure}[tp]
\centering
\includegraphics{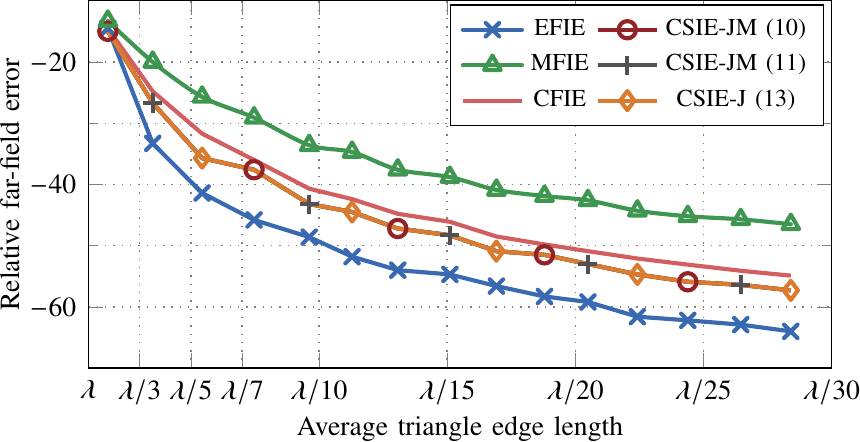}
\caption{Scattered far-field accuracy for the  1\,m square cube.\label{fig:c-acc}}\vspace*{-0.1cm}
\end{figure}

First, the convergence behavior with mesh refinement is analyzed for a 1\,m square cube with plane wave illumination with a wavelength of 2\,m. 
To analyze the MFIE operator within the CSIE on its own, a combination parameter $\alpha=10$ is chosen 
for the mesh refinement analysis. We consider three different formulations of the CSIE: first the straight-forward version CSIE-JM~\eqref{eq:CSIE-saddel} as originally proposed in~\cite{Kornp2017}, then the modified CSIE-JM with weighted magnetic currents according to~\eqref{eq:block}, and finally the novel formulation~\eqref{eq:csiej}. 
In Fig.~\ref{fig:c-conv}, the iterative solver convergence behavior to a residual error of $10^{-4}$ is analyzed. 
The new CSIE-J formulation shows an acceptable convergence, faster than the EFIE and almost as fast as the CFIE. For this small example, in the absence of interior resonances, the CSIE-JM is much slower than the EFIE and, thus, the formulation with the slowest iterative solver convergence in this comparison. 
The block-diagonal preconditioning can only overcome this in part.
In Fig.~\ref{fig:c-acc}, the averaged error of the scattered far-field with respect to a higher-order EFIE solution on the finest mesh is shown. 
As expected, both CSIE-J and CSIE-JM show almost exactly the same error level, which is about 10\,dB better than the error of the MFIE. The EFIE solution is even about 5\,dB more accurate. 
The performance difference of CFIE and CSIE-J will be analyzed in the following more carefully in order to highlight the benefits of the novel CSIE-J. 
It will be shown that only this CSIE is able to outperform the CFIE.

\begin{figure}[!tp]
\centering
\includegraphics{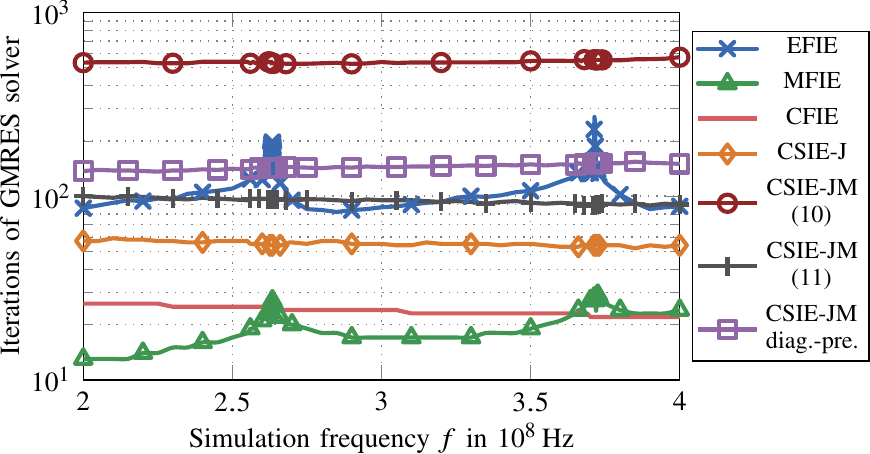}
\caption{Iterative solver convergence analysis for the first two resonances of a PEC sphere with a diameter of 1\,m.\label{fig:s2-conv}}
\end{figure}

\begin{figure}[!tp]
\centering
\includegraphics{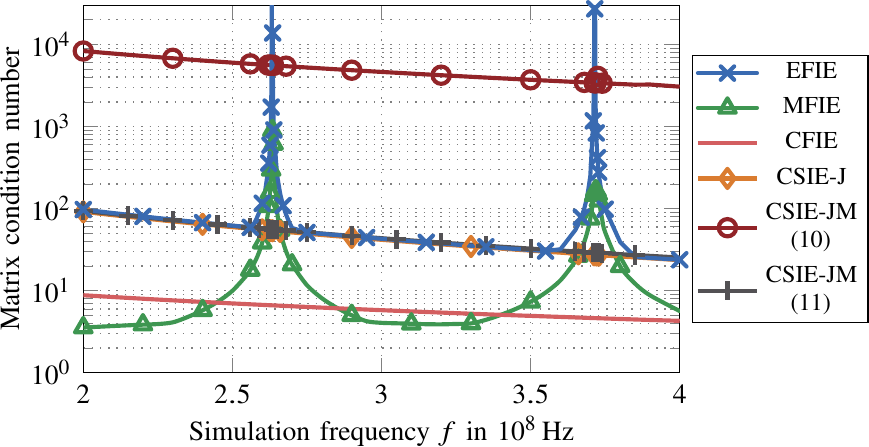}
\caption{Condition number around  the first two resonances of a 1\,m sphere.\label{fig:s2-cond}}
\end{figure}

\begin{figure}[!tp]
\centering
\vspace*{-0.08cm}
\includegraphics{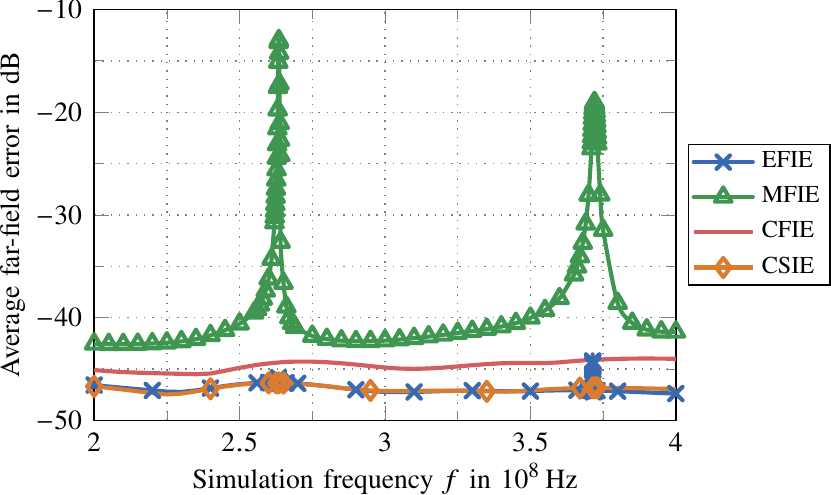}
\caption{Accuracy analysis around the first two resonances of a PEC sphere.\label{fig:s2-acc}}
\end{figure}

A PEC sphere with a diameter of 1\,m and 999 electric current unknowns is considered. 
In the simulated frequency range, the average edge length of the mesh is between $\lambda/25$ and $\lambda/7$. 
The first two interior resonances, corresponding to the resonance frequencies of a spherical PEC cavity with the same diameter, are  observed for the EFIE and the MFIE. 
In Fig.~\ref{fig:s2-conv}, 
the iterative solver convergence of a GMRES solver is analyzed for a residual error of $10^{-4}$. 
In the CSIE-J, the inversion of the Gram matrix for the side condition is solved to a residual error of $10^{-5}$, which is one order of magnitude lower. 
The EFIE, MFIE, CFIE and CSIE-J solutions are obtained without preconditioning; the CSIE-JM~\eqref{eq:CSIE-saddel} is solved without preconditioning, with diagonal preconditioning and with block-diagonal preconditioning according to~\eqref{eq:block}.  
Both EFIE and MFIE suffer from the interior resonances due to a null space at resonance frequencies. 
The CFIE and both CSIE formulations show a stable convergence behavior, while only CSIE-J is faster than the EFIE on average. 
In addition, the matrix condition numbers without preconditioning (except for CSIE-JM) are shown in Fig.~\ref{fig:s2-cond}.  
\begin{figure}[!tp]
\centering
\vspace*{-0.08cm}
\includegraphics{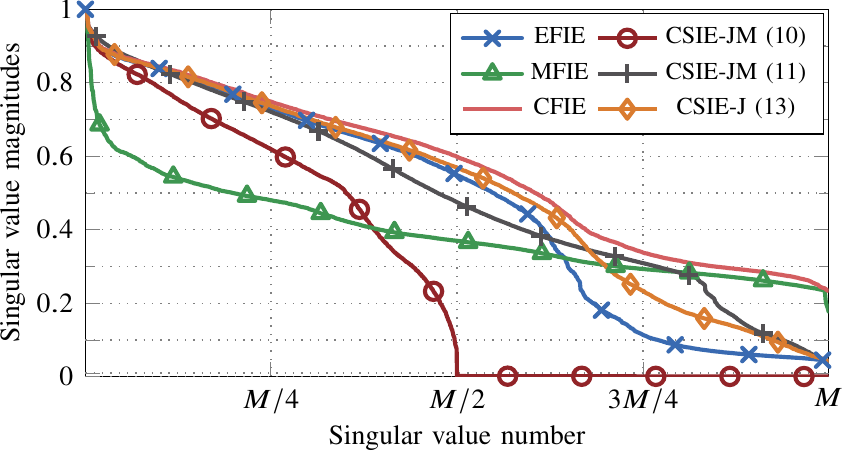}
\caption{Normalized singular value spectra for the PEC sphere at $0.4$\,GHz. For $J$ and $M$ unknowns, $M=1998$, else $M=999$.\label{fig:s2-spec}}
\end{figure}%
Again, the CFIE and CSIE formulations have the advantage of an absolutely stable condition number at interior resonance frequencies. However, one observation is remarkable: 
While the EFIE, CSIE-J and both preconditioned CSIE-JM equations show almost the same condition number, the iterative solver convergence of the CSIE-J is much better. 
The most reasonable cause is a different eigenvalue distribution in the complex plane, e.g. by an improved clustering more suitable for GMRES. 
In the singular-value  spectrum, see Fig.~\ref{fig:s2-spec}, we can recognize a certain amount of singular values with larger values for the CSIE-J as compared to the EFIE. 
The better iterative solver convergence of CSIE-J has been verified for many different right-hand sides, including random vectors and unit vectors of the standard basis. 
Furthermore, the accuracy as compared to the analytical Mie series solution is shown in Fig.~\ref{fig:s2-acc}. 
It is observed that the MFIE shows the poorest accuracy, with huge deteriorations at the resonance frequencies. 
Compared to this behavior, the EFIE only shows inaccuracies quite exactly at the resonance frequencies. 
The CFIE is influenced by the MFIE error and, thus, its accuracy is worse than the CSIE. 
Both CSIE formulations show absolutely the same error level, as expected, which is also about the same as the EFIE error. 

\begin{figure}[!tp]
\centering
\includegraphics[]{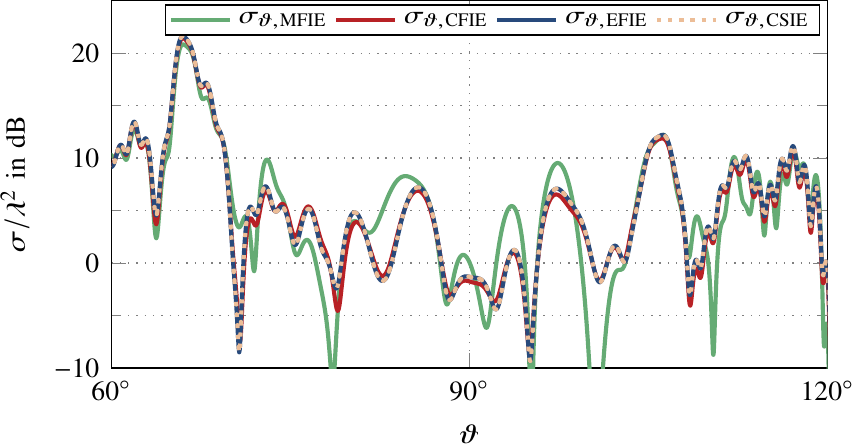}
\caption{Side-lobe section of the bistatic radar cross section~$\sigma$ of the Flamme for the EFIE, CFIE, and CSIE solutions.\label{fig:rcs}}
\end{figure}

As a further problem, the stealth object Flamme is simulated with $\lambda/10$ mean edge length discretization and $287\,922$ electric current unknowns~\cite{GuerelBagciCastelliEtAl2003,Eibert_2005}. 
The MFIE takes $1267$ iterations for convergence to a residual error of $10^{-5}$, the CFIE $315$, the CSIE-J $448$, and the CSIE-JM $1264$ iterations. 
The EFIE did not converge without preconditioning even after thousands of iterations. 
With diagonal preconditioning, the iterative solution of the MFIE is accelerated to $899$ iterations and the iterative solution of the CFIE to $253$. 
The CSIE-J converges within $384$ iterations if only the EFIE operator is taken for preconditioning and within $381$ iterations if the approximated diagonal according to~\eqref{eq:CSIE-diag} is utilized. 
The average time per matrix-vector product is $10.5$ seconds for EFIE, MFIE, and CFIE, and $20.7$ seconds for the CSIE-JM, where all computations have been performed single-threaded on an Intel E5-1650 v4 processor with $3.6$\,GHz clock speed. 
The explicit inversion of the CS condition takes only additional $0.04$ seconds on average, or about $10$ iterations, for a residual error of $5\cdot10^{-7}$. 
Results of the bistatic radar cross section are shown in Fig.~\ref{fig:rcs}. 
Compared to an 3rd order EFIE solution on a refined mesh~\cite{Ismat2009MoM}, the MFIE shows a maximum error of $-37.3$\,dB, the CFIE of $-44.2$\,dB, both CSIE-J and CSIE-JM of $-57.8$\,dB, and the EFIE of $59.8$\,dB.

\begin{figure}[!tp]
\centering
\includegraphics[]{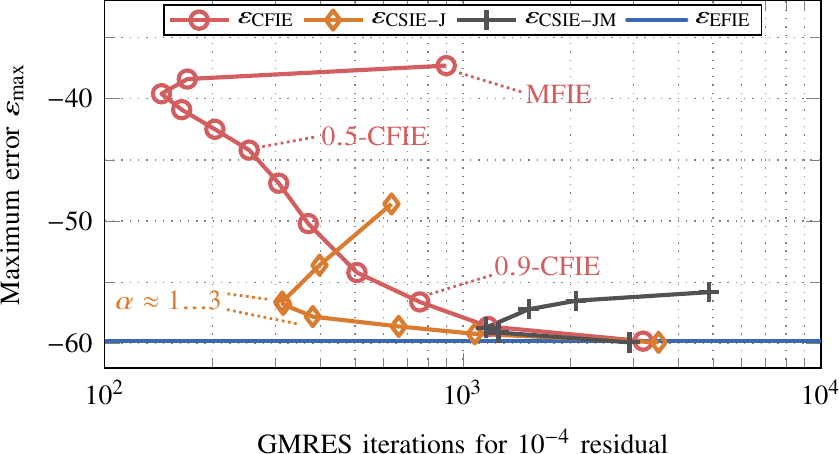}
\caption{Error level vs.\ iterations with different weightings in CFIE and CSIE for Flamme.\label{fig:weight}}
\end{figure}

With a difference in the iterative solver convergence,  it is still not clear whether the CSIE-J or the CFIE offer the better accuracy-convergence trade-off. 
Therefore, we analyze the influence of the corresponding  weighting or combination factors for the cases of CSIE-J and CFIE with diagonal preconditioning and CSIE-JM with block-diagonal preconditioning in Fig.~\ref{fig:weight}. 
The simulated object is again the Flamme with the same parameters as before. 
It becomes clear that the CSIE-JM can provide the same error level only at the cost of more solver iterations.
Additionally, the CSIE-J is able to outperform the CFIE for near-EFIE accuracies. 
For instance, the CSIE-J with $\alpha=2.7$ converges within 313 iterations and offers an error level of  $56.6$\,dB, while the CFIE with combination parameter~$0.9$ (90\% EFIE) converges in 758 iterations for the same error. 
This still means less than half the solver memory consumption for the CSIE-J, approximately the same memory for the matrices, and about twice the time per matrix-vector-product (but less solution time overall).

\section{Conclusion} 

A new form of the CSIE with weak-form CS condition has been introduced, showing an improved iterative solver convergence and a reduced number of unknowns. The magnetic current unknowns are explicitly computed from the electric currents for each matrix-vector product within the iterative solver and, therefore, eliminated. 
The new CSIE is able to offer a better accuracy/convergence trade-off than the CFIE and the weighting factor in the CSIE is not a critical choice (other than in the CFIE) in order to achieve almost EFIE accuracy.

\bibliographystyle{IEEEtran}
\bibliography{IEEEabrv,ref2}
\end{document}